\documentclass[10pt]{amsart}
\usepackage{amssymb, amsmath}
\usepackage[mathscr]{eucal}
\usepackage{graphicx}

\begin{document}

\title[Flexible polyhedra: algebra vs analysis]{Algebra versus analysis in the theory\\ of flexible polyhedra}
\author[Victor Alexandrov]{Victor Alexandrov}
\address{
Sobolev Institute of Mathematics, Koptyug ave., 4, Novosibirsk, 630090, Russia and
Department of Physics, Novosibirsk State University, Pirogov str., 2, Novosibirsk, 630090, Russia}
\email{alex@math.nsc.ru}
\thanks{The author was supported by the Federal Program `Research and educational 
resourses of innovative Russia in 2009--2013' (contract~02.740.11.0457) and the
Russian Foundation for Basic Research (grant~10--01--91000--ANF)}
\date{July 6, 2010}
\begin{abstract}
Two basic theorems of the theory of flexible polyhedra were proven by completely different methods:
R.~Alexander used analysis, namely, the Stokes theorem, to prove that the total mean curvature remains
constant during the flex, while I.Kh.~Sabitov used algebra, namely, the theory of resultants,
to prove that the oriented volume remains constant during the flex.
We show that none of these methods can be used to prove the both theorems.
As a by-product, we prove that the total mean curvature of any polyhedron in the Euclidean 3-space
is not an algebraic function of its edge lengths. 
\end{abstract}
\maketitle

A polyhedron (more precisely, a polyhedral surface)
is said to be \textit{flexible} if its spatial shape can 
be changed continuously due to changes of its dihedral angles only,
i.\,e., if every face remains congruent to itself during the flex.
In other words, a polyhedron $P_0$ is flexible if it is included
in a continuous family $\{P_t\}$, $0\leqslant t\leqslant 1$, of
polyhedra $P_t$ such that, for every $t$, the corresponding faces of $P_0$ and $P_t$
are congruent while the polyhedra $P_0$ and $P_t$ are not congruent.
The family $\{P_t\}$, $0\leqslant t\leqslant 1$, is called the flex of $P_0$.
Self-intersections are possible both for $P_0$ and $P_t$  
provided the converse is not formulated explicitly.
Without loss of generality we assume that the faces of the
polyhedra are triangular.

Flexible self-intersection free sphere-homeomorphic polyhedra
in Euclidean 3-space were constructed by R.~Connelly
thirty years ago \cite{Co80}, \cite{Ku79}.
Since that time, various non-trivial properties of flexible 
polyhedra were discovered in the Euclidean 3-space \cite{Sc04}
and 4-space \cite{St00} (for results in the hyperbolic 3-space, see \cite{St06}).
Let us formulate two of them in a form suitable for our purposes.

Let $P$ be a closed oriented polyhedron in $\mathbb R^3$, let
$E$ be the set of its edges, let $|\ell |$ be the length of the edge $\ell$,
and let $\alpha (\ell)$ be the dihedral angle of $P$ at the edge $\ell$
measured from inside of $P$.
The sum
$$
M(P)=\frac12\sum_{\ell\in E} | \ell | \bigl(\pi -\alpha (\ell)\bigr) 
$$
is called the \textit{total mean curvature} of $P$.

\textbf{Theorem 1.}
\textit{Let $P_0$ be a flexible polyhedron in $\mathbb R^3$ and let $\{P_t\}$, $0\leqslant t\leqslant 1$, be its flex.
The total mean curvature $M(P_t)$ is independent of $t$.}\hfill$\square$

Theorem 1 was proved by R.~Alexander \cite{Al85}
for all Euclidean $n$-spaces, $n\geqslant 3$, though no example of 
a flexible polyhedron is known for $n\geqslant 5$.

Let $P_0$ be a flexible polyhedron and let $\{P_t\}$, $0\leqslant t\leqslant 1$, be its flex.
Let $\boldsymbol{r}_t$ be the point of the polyhedron $P_t$ 
which corresponds to the point $\boldsymbol{r}_0\in P_0$. 
It follows from the definition of a flexible polyhedron that, for any curve 
$\boldsymbol{\gamma}_0\subset P_0$, the length of the curve 
$\boldsymbol{\gamma}_t=\{\boldsymbol{r}_t\vert\boldsymbol{r}_0\in\boldsymbol{\gamma}_0\}\subset P_t$
is independent of $t$. 
The reader can easily verify the well-known fact that, for any curve $\boldsymbol{\gamma}_0\subset P_0$,
the length of the curve 
$\widetilde{\boldsymbol{\gamma}}_t=\{\boldsymbol{r}_0+t\boldsymbol{v}\vert \boldsymbol{r}_0\in\boldsymbol{\gamma}_0\}$
is stationary at $t=0$, where $\boldsymbol{v}=\frac{d}{dt}\vert_{{}_{\scriptstyle{t=0}}}\boldsymbol{r}_t$
is the velocity vector of the point $\boldsymbol{r}_t$ at $t=0$. 
Obviously, the vector field $\boldsymbol{v}$ is linear on every face of $P_0$.
This leads to the following well-known definition: 
a vector field $\boldsymbol{w}$ on a polyhedron $P$, which  is linear on every face of $P$, 
is said to be its \textit{infinitesimal flex} if, for any curve $\boldsymbol{\gamma}\subset P$,
the length of the curve 
$\widetilde{\boldsymbol{\gamma}}(t)=\{\boldsymbol{r}+t\boldsymbol{w} \vert\boldsymbol{r}\in\boldsymbol{\gamma}\}$
is stationary at $t=0$, see \cite{Sa92} for more detail.
Of course, the velocity vector field of a flexible polyhedron is its infinitesimal flex 
(but the converse is not necessarily true).

In \cite{Al85} Theorem 1 was obtained as an obvious corollary of the following Theorem 2
proved for $\mathbb R^n$, $n\geqslant 3$. 

\textbf{Theorem 2.}
\textit{Let $P$ be a closed oriented polyhedron in $\mathbb R^3$,
let $\boldsymbol{w}$ be its infinitesimal flex, and
let $P(t)=\{\boldsymbol{r}+t\boldsymbol{w} \vert\boldsymbol{r}\in P\}$.
Then}
$\frac{d}{dt}\vert_{{}_{\scriptstyle{t=0}}}M\bigl(P(t)\bigr)=0.$\hfill$\square$

In \cite{Al85} Theorem 2 was proved with the help of the Stokes theorem.
Later several authors observed (see, e.\,g., \cite{AL98}, \cite{So04})
that Theorem 2 follows immediately from the Schl\"afli differential formula, 
which, in turn, is based on the Stokes theorem.
In any case, all known proofs of Theorems 1 and 2 belong to Analysis.

In \cite{Sa96} I.Kh.~Sabitov proved another highly non-trivial property
of the flexible polyhedra that may be formulated as follows and
whose many-dimensional analog is not known yet.

\textbf{Theorem 3.}
\textit{If $\{P_t\}$ is a flex of an orientable polyhedron in $\mathbb R^3$,
then the oriented volume of $P_t$ is constant in $t$.}\hfill$\square$

In \cite{Sa96}, \cite{Sa98}, and later in \cite{CSW7} Theorem 3
was obtained as an obvious corollary of the following Theorem 4
valid since every continuous mapping, whose image is a finite set,
is constant.

\textbf{Theorem 4.}
\textit{For the set $\mathscr P_K$ of all (not necessarily flexible)
closed polyhedra in $\mathbb R^3$ with triangular faces and with a prescribed 
combinatorial structure $K$ there exists a universal polynomial $\mathfrak{p}_K$ 
of a single variable whose coefficients are universal polynomials
in the edge lengths of a polyhedron $P\in\mathscr P_K$ and
such that the oriented the volume of any $P\in\mathscr P_K$
is a root of $\mathfrak{p}_K$.}\hfill$\square$

In \cite{Sa96} and \cite{Sa98} Theorem 4 was proved with the help of the theory
of resultants, while in \cite{CSW7} it was proved with the help of the theory
of places. In any case, all known proofs of Theorems 3 and 4 belong to Algebra.

Now recall that
the derivative of the volume of a deformable domain in $\mathbb R^3$ equals
one third of the flux across the boundary of the velocity vector of a point 
of the boundary of the domain, i.\,e., equals one third of the integral 
of the normal component of the velocity over the boundary.
Hence, an `infinitesimal  version' of Theorem~3 should read 
that, \textit{for every orientable polyhedron $P$, the flux
$
\int_P (\boldsymbol{w},\boldsymbol{n})\, dP
$
of any infinitesimal flex $\boldsymbol{w}$ equals zero.}
Theorem~5 shows that this is not the case and, thus, 
that Theorem~3 can not be proved by means of Analysis 
like Theorem~1.

\textbf{Theorem 5.}
\textit{There is a closed oriented polyhedron~$P$ 
in~$\mathbb R^3$ with the following properties}: 

(i)~\textit{the flux across $P$ of some infinitesimal flex 
$\boldsymbol{w}$ of $P$ is non-zero};

(ii)~\textit{$P$ contains no vertex $V$ whose star lies in 
a plane}; 

(iii)~\textit{$P$ contains no vertex $V$ such that some three edges
of $P$ incident to $V$ lie in a plane.}

\textit{Proof.}
Let $T_0=ABCD$ be an arbitrary tetrahedron in $\mathbb R^3$.
Let a polyhedron $T_1$ be obtained from $T_0$ by triangulation of 
the face $ABC$ that uses one additional vertex $V$ (see Fig.~1). 
Let a vector field $\boldsymbol{w}_1$ be linear on each face
of $T_1$, be equal zero at the vertices $A$, $B$, $C$, and $D$,
and be equal a non-zero vector perpendicular to the face~$ABC$ 
at the vertex~$V$. It is easy to check that $\boldsymbol{w}_1$
is an infinitesimal flex of $T_1$ and its flux across $T_1$
is non-zero (because the scalar product   
$(\boldsymbol{w}_1,\boldsymbol{n})$ 
is either everywhere non-negative or everywhere non-positive;
here $\boldsymbol{n}$ stands for the unit normal vector field 
on $T_1$). 
Hence, $T_1$ satisfies the condition~(i), though it does not 
satisfy the conditions~(ii) and~(iii).
\begin{figure}
\includegraphics[width=0.4\textwidth]{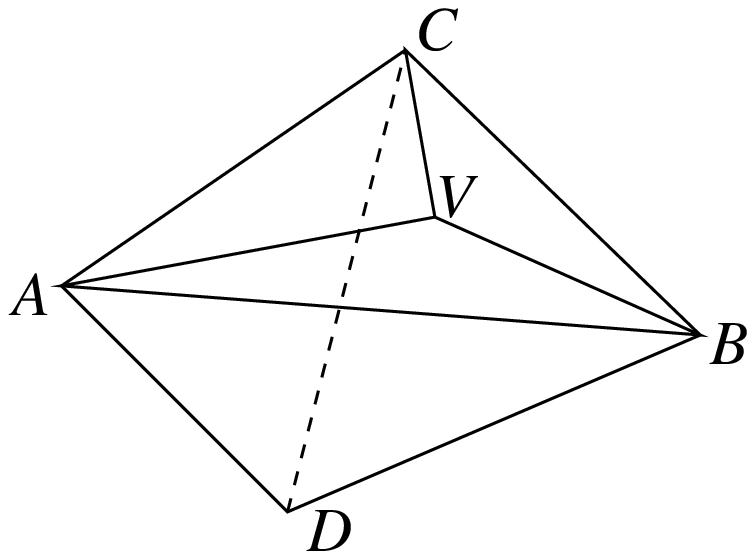}\hfill
\includegraphics[width=0.4\textwidth]{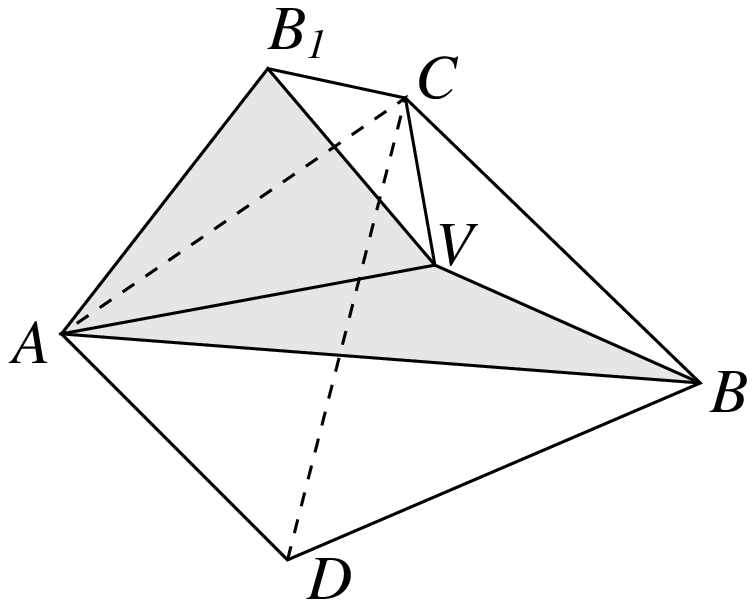}\\
\hskip-65mm\parbox[t]{0.47\textwidth}{\caption{Polyhedron~$T_1$}}\hskip10mm
\parbox[t]{0.47\textwidth}{\caption{Polyhedron~$T_2$}}
\end{figure}

Let $B_1$ be an arbitrary point in $\mathbb R^3$
which does not lie in the plane $ABC$.
Replace the face $ACV$ of the polyhedron $T_1$ with
the lateral surface of the triangular pyramid $ACVB_1$
with the base $ACV$.
Denote the resulting polyhedron by $T_2$ (see Fig.~2).
It is easy to check that the infinitesimal flex $\boldsymbol{w}_1$
of $T_1$ can be extended in a unique way to an infinitesimal flex
of $T_2$. Denote it by $\boldsymbol{w}_2$ and observe that the flux
of~$\boldsymbol{w}_2$ across~$T_2$ equals the flux
of~$\boldsymbol{w}_1$ across~$T_1$ (because the flux of
$\boldsymbol{w}_2$ across the triangular pyramid~$ACVB_1$ 
equals zero). Hence,~$T_2$ satisfies the conditions~(i)
and~(ii), though it does not satisfy the condition (iii).

Recall the construction of the Bricard flexible octahedron 
of type 1.
Consider a disk-homeomorphic piece-wise linear surface $S$ 
in $\mathbb R^3$  composed of four triangles 
$ABV$, $BA_1V$, $A_1B_1V$, and $B_1AV$ 
such that
$|AB|=|A_1B_1|$ and $|B_1A|=|BA_1|$ (see Fig.~3). 
It is known~\cite{Ku79}  that such a spatial quadrilateral
$ABA_1B_1$ is symmetric with respect to the line~$L$ passing 
through the middle points of its diagonals 
$AA_1$ and $BB_1$.
Glue together~$S$ and its symmetric image with respect to~$L$ 
(see Fig.~4).
The resulting polyhedral surface with 
self-intersections is flexible (because $S$ is flexible) 
and is combinatorially equivalent to the surface of the 
regular octahedron. 
This is known as the Bricard octahedron of type 1.
\begin{figure}
\includegraphics[width=0.49\textwidth]{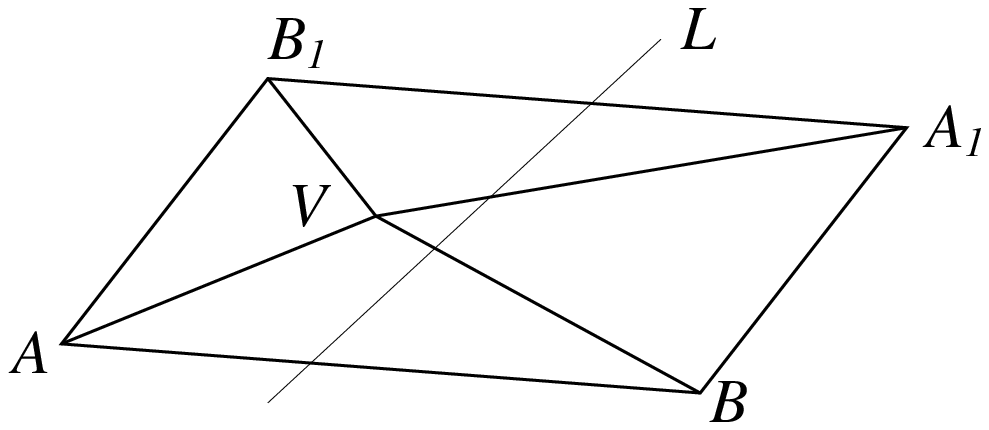}\hfill
\includegraphics[width=0.49\textwidth]{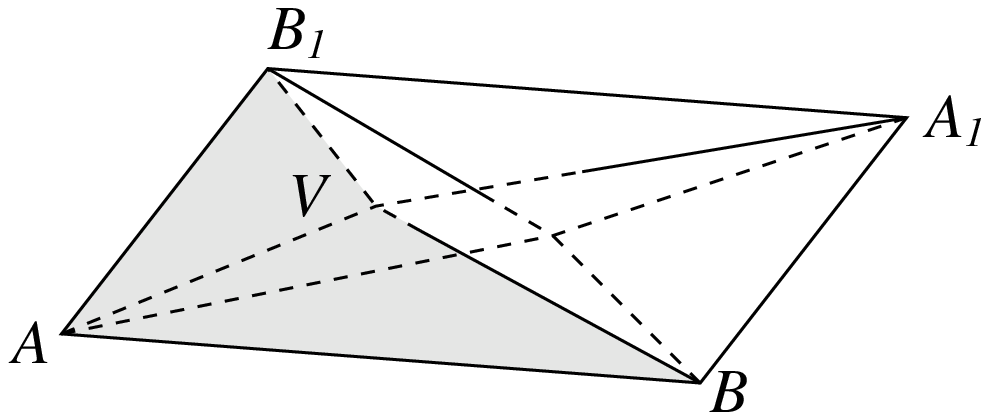}\\
\hskip-75mm\parbox[t]{0.47\textwidth}{\caption{Disk~$S$}}\hskip20mm
\parbox[t]{0.47\textwidth}{\caption{Bricard octahedron}}
\end{figure}

Construct the Bricard octahedron of type 1 such that
its triangles $ABV$ and $AVB_1$ are congruent with the
corresponding triangles  $ABV$ and $AVB_1$ of the 
polyhedron~$T_2$.
Remove those triangles from $T_2$ and from the Bricard 
octahedron and glue together the remaining parts of 
those polyhedra.
Observe that the resulting polyhedron $P$ satisfies 
the conditions~(i)--(iii).
In fact, the infinitesimal flex $\boldsymbol{w}_2$
of $T_2$ can be extended in a unique way to an 
infinitesimal flex of $P$. 
Denote it by $\boldsymbol{w}$ and observe that the flux
of  $\boldsymbol{w}$ across~$P$ equals the flux of
 $\boldsymbol{w}_2$ across~$T_2$ (because 
the restriction of the flux $\boldsymbol{w}$ 
on the Bricard octahedron is generated by some its flex
and, thus equals zero, since, e.\,g., the oriented volume
of the Bricard octahedron is known to be constant
for every flex \cite{CSW7}, \cite{Sa92}--\cite{Sc04}).
\hfill$\square$

Observe that an `algebraic version' of Theorem~2 should read
that \textit{the total mean curvature of any oriented 
polyhedron is a root of some universal polynomial $\mathfrak{p}$
in a single variable whose coefficients are universal polynomials
in the edge lengths of the polyhedron.} 
This statement would, obviously, imply Theorem~1.
But this statement also implies that the total mean curvature
of a polyhedron is an algebraic function of its edge lengths.
Theorem~6 shows that this is not the case and, thus, 
that Theorem~1 can not be proved by means of Algebra 
like Theorem~3.

\textbf{Theorem 6.}
\textit{The total mean curvature of any closed oriented 
polyhedron in~$\mathbb R^3$ is not an algebraic function
of its edge lengths.}

\textit{Proof.} 
For every $K$, the class of
all polyhedra of combinatorial type $K$
contains a family $\varDelta_K$ of polyhedra
depending on a single independent variable $l>0$
which may be described as follows (see Fig.~5):
there are two triangles $ABC$ and $ACD$ in $\mathbb R^3$
such that 

($\alpha$) $|BC|=|CD|=|BD|=1$ and 
$|AB|=|AC|=|AD|=l$; 

($\beta$) every $P\in\varDelta_K$ has the triangles $ABC$ 
and $ACD$  as their faces; 

($\gamma$) every vertex of $P\in\varDelta_K$  either
coincides with $A$, $B$, $C$, or $D$
or is an interior point of the segment $BD$
or is an interior point of the triangle $BCD$;

($\delta$) the set of the interior points of every edge of 
$P\in\varDelta_K$  either coincides with $AB$, $AC$, or $AD$
or lies in the open segment $BD$ or lies in
one of the open triangles $ABD$ or $BCD$.

A direct calculation shows that, for every $P\in\varDelta_K$,
the total mean curvature $M(P)$ of $P$ equals
\begin{align*}
M(P)=&\frac32\bigl(\pi-\varphi(l)\bigr)+
\frac32 l\bigl(\pi-\psi(l)\bigr)\\
=&\frac32\biggl(\pi-\arccos\frac{1}{2\sqrt{3}\sqrt{4l^2-1}}\biggr)+
\frac32 l\biggl(\pi-\arccos\frac{2l^2-1}{4l^2-1}\biggr),
\end{align*}
where $\varphi(l)$ is the inner dihedral angle 
of $P$ at the edge $BC$ and 
$\psi(l)$ is the inner dihedral angle 
of $P$ at the edge $AC$.
Consider the right-hand side of the last formula 
as a function of complex variable~$l$.
Obviously, this function has a non-algebraic singularity
(known also as a logarithmic branch point) over
$l=0$ (as well as over $l=\pm1/\sqrt{3}$
and $l=\pm i \sqrt{48/13}$).
\begin{figure}
\includegraphics[width=0.4\textwidth]{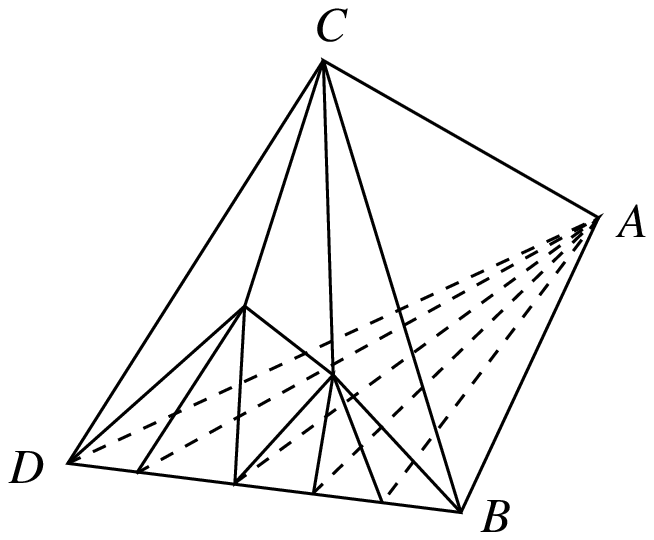}
{\caption{Polyhedron of the class~$\varDelta_K$}}
\end{figure}

Recall that a function $w=f(z)$ of a single complex 
variable $z$ is called algebraic, if 
there is a polynomial $\mathfrak{p}(w,z)$ in two 
variables which does not vanish identically
and such that $\mathfrak{p}(f(z),z)\equiv 0$. 
It is known that an analytic function
of a single complex variable is an algebraic function 
if and only if it has a finite number of branches and 
at most 
algebraic singularities \cite[p.~306]{Ah79}.
Hence, the total mean curvature~$M(P)$
is not an algebraic function 
of the variable $l$ for the polynomials~$P$
of the class $\varDelta_K$.
Thus,~$M(P)$ is not an algebraic function
of the edge lengths of the polyhedron $P$.  
\hfill$\square$

\end{document}